\documentclass[12pt]{amsart}
\usepackage{amsmath,amssymb,amscd,graphicx}
\usepackage[all]{xy}
\usepackage[latin1]{inputenc}
\usepackage[T1]{fontenc}
\usepackage{hyperref}
\usepackage{ifthen}
\usepackage{fullpage}
\usepackage{subsupscripts}
\usepackage{xcolor}

\newboolean{workmode}
\setboolean{workmode}{false}

\newboolean{sections}
\setboolean{sections}{true}

\input xy
\xyoption{all}

\def\eoe{\unskip\ \hglue0mm\hfill$\diamond$\smallskip\goodbreak}

\makeatletter
\def\th@plain{%
  \thm@notefont{}
  \itshape 
}
\def\th@definition{%
  \thm@notefont{}
  \normalfont 
}
\makeatother

\newcommand{\calD}{{\mathcal{D}}}

\newcommand{\calO}{{\mathcal{O}}}

\newcommand{\RR}{{\mathbb{R}}}  

\newcommand{\Ad}{{\operatorname{Ad}}}  
\newcommand{\Bi}{{\mathbf{Bi}}} 
\newcommand{\Diffeol}{\mathbf{Diffeol}} 
\newcommand{\pr}{{\operatorname{pr}}} 

\newcommand{\toto}{{~\rightrightarrows~}} 
\newcommand{\hook}{{\lrcorner\,}} 

\newcommand{\ifwork}[1]{\ifthenelse{\boolean{workmode}}{#1}{}}
\newcommand{\comment}[1]{}
\newcommand{\mute}[1]{}
\newcommand{\printname}[1]{}

\ifwork{
\renewcommand{\comment}[1]{{\marginpar{*}\ \scriptsize{#1}\ }}
\renewcommand{\printname}[1]
    {\smash{\makebox[0pt]{\hspace{-1.0in}\raisebox{8pt}{\tiny #1}}}}
}

\newcommand{\labell}[1] {\label{#1} \printname{#1}}

\newcommand{\ifsection}[2]{\ifthenelse{\boolean{sections}}{#1}{#2}}

\theoremstyle{plain}
\ifsection{
    \newtheorem{theorem}{Theorem}[section]
}
{
    \newtheorem{theorem}{Theorem}
}

\newtheorem{proposition}[theorem]{Proposition}

\newtheorem{corollary}[theorem]{Corollary}
\newtheorem{lemma}[theorem]{Lemma}

\theoremstyle{definition}
\newtheorem{definition}[theorem]{Definition}
\newtheorem{example}[theorem]{Example}
\newtheorem{remark}[theorem]{Remark}

\definecolor{jaw}{rgb}{0,.5,0}

\newcommand{\basic}{{\operatorname{basic}}} 

\author{Jordan Watts}

\address{Department of Mathematics\\ Central Michigan University\\ Mount Pleasant\\ MI 48859}
\email{jordan.watts@cmich.edu}

\title{The Orbit Space and Basic Forms of a Proper Lie Groupoid}
\date{\today}

\begin{document}

\maketitle

\begin{abstract}
A classical result in differential geometry states that for a free and proper Lie group action, the quotient map to the orbit space induces an isomorphism between the de Rham complex of differential forms on the orbit space and the basic differential forms on the original manifold.  In this paper, this result is generalized to the case of a proper Lie groupoid, in which the orbit space is equipped with the quotient diffeological structure.  As an application of this, we obtain a de Rham theorem for the de Rham complex on the orbit space.
\end{abstract}

{\small\emph{2010 Mathematics Subject Classification.}  Primary 58H05, Secondary 22A22\\
\indent Keywords: diffeology, Lie groupoid, de Rham complex, basic forms, linearization}

\section{Introduction}\labell{s:intro}

Given a proper Lie group action on a smooth manifold, the Slice Theorem of Koszul and Palais \cite{koszul,palais} states that any orbit of the action has an invariant tubular neighbourhood equivariantly diffeomorphic to a linear model about the orbit.  From this it follows that in the case of a free and proper action, the quotient map induces an isomorphism between the complex of differential forms on the orbit space and the complex of basic differential forms on the manifold.

Using the theory of diffeology, which generalises the theory of smooth manifolds (see Definition~\ref{d:diffeology}), this result was extended by the author to include any compact Lie group action (not necessarily free) \cite[Chapter 3]{watts-phd}, and further to Lie group actions for which the identity component of the Lie group acts properly by Karshon-Watts \cite{KW}.  The purpose of this note is to extend the result to proper Lie groupoids; namely, to prove Theorem~\ref{t:main} below.  Denote by $\Bi$ the bicategory of Lie groupoids with (right principal) bibundles as arrows and isomorphisms of bibundles as $2$-arrows (see \cite{lerman,MM05} for definitions and more details).  Given a Lie groupoid $G=(G_1\toto G_0)$, we define a differential form $\alpha$ on $G_0$ to be basic if $s^*\alpha=t^*\alpha$, where $s$ and $t$ are the source and target maps. 

\begin{theorem}\labell{t:main}
Given a proper Lie groupoid $G=(G_1\toto G_0)$, the quotient map $G_0\to G_0/G_1$ induces an isomorphism between the diffeological de Rham complex of the orbit space $G_0/G_1$ and the complex of basic forms on $G_0$; moreover, this isomorphism is natural on the full sub-bicategory of proper Lie groupoids in $\Bi$.
\end{theorem}

The natural isomorphism in Theorem~\ref{t:main} suggests that functors should be involved.  Indeed, we show that there is a $2$-functor $\Psi$ from $\Bi$ to the category of diffeological spaces, sending a Lie groupoid $G$ to its orbit space $G_0/G_1$, a bibundle $P\colon G\to H$ to a smooth map $\Psi_P\colon G_0/G_1\to H_0/H_1$, and a $2$-arrow between bibundles to a trivial $2$-arrow.  If $\Omega^*_\basic$ is the (contravariant) functor sending a Lie groupoid to its complex of basic forms, and $\Omega^*$ is the functor sending a diffeological space to its de Rham complex, then Theorem~\ref{t:main} states that when restricting to proper groupoids, there is a natural isomorphism from $\Omega^*\circ \Psi$ to $\Omega^*_\basic$ given by the pullback map induced by the quotient map from the base of a groupoid to its orbit space.

The proof of Theorem~\ref{t:main} relies on the corresponding result for compact Lie group actions (see Theorem~\ref{t:KW}), and a linearisation theorem for proper Lie groupoids (see Theorem~\ref{t:linearisable}).  This linearisation theorem has been developed throughout a series of works, which include authors such as Zung, Weinstein, Crainic-Struchiner, and del Hoyo-Fernandes \cite{zung,weinstein1,weinstein2,CS,dHF}; for our purposes we adopt the language of del Hoyo-Fernandes.

We obtain some immediate consequences of Theorem~\ref{t:main}.  It follows from the Slice Theorem for a proper Lie group action on a manifold that the cohomology of the basic forms is isomorphic to the singular cohomology on the orbit space.  This result was extended to proper Lie groupoids by Pflaum-Posthuma-Tang \cite[Section 8]{PPT}.  Thus together with Theorem~\ref{t:main}, we obtain a de Rham theorem for the orbit space of a proper Lie groupoid which is \emph{intrinsic}, in the sense that it only depends on the diffeology (\emph{i.e.}\ smooth structure) of the orbit space, and not on the original Lie groupoid.

\begin{corollary}\labell{c:main}
Given a proper Lie groupoid $G=(G_1\toto G_0)$, the de Rham cohomology of the orbit space $G_0/G_1$ is isomorphic to the singular cohomology of $G_0/G_1$.
\end{corollary}

Another application is a reinterpretation of Corollary~\ref{c:basic forms of the relation groupoid} and Theorem~\ref{t:main} in terms of exactness in the following sequence:
\begin{equation}\labell{e:intro}
\xymatrix{
0 \ar[r] & \Omega^*(G_0/G_1) \ar[r]^{\pi_G^*} & \Omega^*(G_0) \ar[r]^{s^*-t^*} & \Omega^*(G_1)
}
\end{equation}

\begin{corollary}\labell{c:exactness}
For any Lie groupoid $G=(G_0\toto G_1)$, the sequence \ref{e:intro} is exact at $\Omega^*(G_0/G_1)$.  If $G$ is proper, then the sequence is also exact at $\Omega^*(G_0)$.
\end{corollary}

At this point it is natural for the reader to ask whether the properness condition on a Lie groupoid in Theorem~\ref{t:main} can be relaxed.  At this time, the author is not aware of an example in which the isomorphism does \emph{not} hold.  Indeed, even when the quotient has trivial topology but non-trivial diffeology, such as the $1$-dimensional irrational torus, Theorem~\ref{t:main} still holds; see \cite[Example 5.11]{KW}.

This paper is broken down as follows: Section~\ref{s:background} reviews the theory of diffeology.  Section~\ref{s:basic differential forms} reviews and develops needed results on basic differential forms on Lie groupoids from a categorical perspective.  Section~\ref{s:linearisation} reviews linearizations, and proves Theorem~\ref{t:main}.

We end this introduction with a brief survey of related ideas in the literature.  For instance, the functor $\Psi$ between Lie groupoids and diffeological spaces appears in the literature already, sometimes in disguise.  It follows from \cite{IZKZ} that the restriction of $\Psi$ to effective orbifolds, viewed as effective proper \'etale Lie groupoids, is essentially injective on objects.  There is a similar functor to (Sikorski) differential spaces that factors through $\Psi$; see \cite[Theorem B]{watts-orb} for details.  The point-of-view of stacks over manifolds is taken in \cite{WW}, in which a functor sending a stack to its underlying diffeological coarse moduli space is constructed; as differentiable stacks are represented by Lie groupoids, this is just a stacky manifestation of the functor $\Psi$ above extended to all stacks.  Finally, a more detailed study of when and how a diffeologically smooth map between orbit spaces is in the image of the functor $\Psi$ appears in \cite{KZ}.

The author would like to thank Rui Loja Fernandes, Eugene Lerman, and Ioan M\u{a}rcu\c{t} for many illuminating discussions about Lie groupoids, and Yael Karshon for her comments and encouragement to publish this note.

\section{Background}\labell{s:background}

It is assumed that the reader is familiar with Lie groupoids and (right principal) bibundles between them; see \cite{lerman,MM05} for an exposition on these.  For the purposes of this paper, all Lie groupoids $G=(G_1\toto G_0)$ are assumed to be finite dimensional, paracompact, and Hausdorff. Groupoid homomorphisms $F\colon G\to H$ are denoted using $F_0\colon G_0\to H_0$ for the map between bases, and $F_1\colon G_1\to H_1$ for the map between arrow spaces.  A more detailed exposition of diffeological spaces appears in \cite{iglesias}, but we give a brief review of them now.

\begin{definition}[Diffeology] \labell{d:diffeology}
Let $X$ be a set.  A \textbf{parametrisation} of $X$ is a
map of sets $p \colon U \to X$ where $U$ is an open subset of some Euclidean space (no fixed dimension).  A \textbf{diffeology} $\calD$ on $X$ is a family of parametrisations satisfying:
\begin{enumerate}
\item \textbf{(Covering)} $\calD$ contains all constant maps into $X$.
\item \textbf{(Locality)} Let $p\colon U\to X$ be a parametrisation, $\{U_\alpha\}$ an open cover of $U$, and $\{p_\alpha\colon U_\alpha\to X\}\subseteq\calD$ such that $p|_{U_\alpha}=p_\alpha$ for each $\alpha$.  Then $p\in\mathcal{D}$.
\item \textbf{(Smooth Compatibility)} For any $(p\colon U\to X)\in\calD$ and any smooth parametrisation $f\colon V\to U$ of $U$, the composition $p\circ f$ is in $\calD$.
\end{enumerate}
A set $X$ equipped with a diffeology $\mathcal{D}$ is called a
\textbf{diffeological space}, and is denoted by $(X,\mathcal{D})$.
When the diffeology is understood, we will drop the symbol $\mathcal{D}$.
The parametrisations $p\in\mathcal{D}$ are called \textbf{plots}.
Given two diffeological spaces $(X,\mathcal{D}_X)$ and $(Y,\mathcal{D}_Y)$, a map $\varphi \colon X \to Y$ is
\textbf{(diffeologically) smooth} if $\varphi \circ p \in \mathcal{D}_Y$ for any plot $p \in \mathcal{D}_X$.
\end{definition}

Diffeological spaces with smooth maps between them form a category $\Diffeol$, which contains the category of smooth manifolds as a full subcategory. Indeed, the  \textbf{standard diffeology} on $M$ is the set of all smooth parametrisations $f\colon U\to M$, where ``smooth'' is taken in the usual sense.

The usefulness of $\Diffeol$ is illustrated by the fact that it is complete and co-complete \cite{BH,iglesias}.  For instance, given a diffeological space $X$ and an equivalence relation $\sim$ on $X$ with quotient map $\pi\colon X\to X/\!\sim$, the quotient comes equipped with the \textbf{quotient diffeology}.  This is the family of all parametrisations $p\colon U\to X/\!\sim$ satisfying the condition that each point $u\in U$ has an open neighbourhood $V$ such that $p|_V$ factors through $\pi$.  As another example, given a subset $Z\subseteq X$, the \textbf{subset diffeology} on $Z$ consists of all plots of $X$ whose image lies in $X$.  Finally, given another diffeological space $Y$, the \textbf{product diffeology} on $X\times Y$ is the collection of all parametrisations $(p_1,p_2)$ such that $p_1$ is a plot of $X$ and $p_2$ is a plot of $Y$.

\begin{example}\labell{x:orbitspace}
Let $G$ be a Lie groupoid and fix $x\in G_0$. The \textbf{orbit} $\mathcal{O}$ of $G$ through $x$ is the set $$\mathcal{O}=\{y\in G_0\mid\exists g\in G_1 \text{ so that } s(g)=x \text{ and } t(g)=y\}.$$  The subset diffeology on $\mathcal{O}$ gives it the structure of an immersed submanifold of $G_0$ (see \cite[Section 1.2]{CS}).  The \textbf{orbit space} $G_0/G_1$ of $G$ is the quotient of $G_0$ by the equivalence relation $\sim$ given by: $x\sim y$ if $x$ and $y$ are in the same orbit.  $G_0/G_1$ comes equipped with the quotient diffeology induced by the standard manifold diffeology on $G_0$; denote the quotient map by $\pi_G\colon G_0\to G_0/G_1$.

The fibred product $G_0\times_\pi G_0$ is the arrow space of the \textbf{relation groupoid} of $\sim$ $$\{(x_1,x_2)\in G_0\times G_0~|~\exists g\in G_1 \text{ such that } s(g)=x_1\text{ and } t(g)=x_2\}$$ whose source and target are restrictions of the first and second projection maps, respectively, to $G_0$. The arrow space comes equipped with the subset diffeology induced by $G_0\times G_0$, with respect to which the map $(s,t)\colon G_1\to G_0\times_\pi G_0$ is a smooth surjection.
\eoe
\end{example}

\begin{definition}[Differential Forms]\labell{d:forms}
A \textbf{differential $k$-form} $\alpha$ on a diffeological space $X$ is an assignment to each plot $p\colon U\to X$ a differential form $\alpha_p\in\Omega^k(U)$, satisfying $\alpha_{p\circ f}=f^*\alpha_p$ for any smooth parametrisation $f\colon V\to U$ of $U$.
Denote the collection of $k$-forms on $X$ by $\Omega^k(X)$.  The \textbf{exterior derivative} $d\colon\Omega^k(X)\to\Omega^{k+1}(X)$ is defined plot-wise: $(d\alpha)_p=d(\alpha_p).$
\end{definition}

The $0$-forms of a diffeological space are exactly the smooth functions $f\colon X\to\RR$; the collection of all differential forms on $X$ forms a de Rham complex $(\Omega^*(X),d)$, although we will abbreviate the notation to just $\Omega^*(X)$; and a smooth map $\varphi\colon X\to Y$ induces a \textbf{pullback map} $f^*\colon\alpha\mapsto f^*\alpha$, which is a map of complexes.  In particular, we have a contravariant functor $\Omega$ sending a diffeological space $X$ to the complex $\Omega^*(X)$ and a smooth map to its corresponding pullback map.  Finally, on a manifold equipped with the standard diffeology, the diffeological de Rham complex is exactly the standard one.

\section{Basic Differential Forms and Bibundles}\labell{s:basic differential forms}

In this section, we collect a number of results regarding basic forms of a Lie groupoid.  After defining them and relating them to the classical notion for a Lie group action (Lemma~\ref{l:basic forms}), we prove that pullbacks of forms from an orbit space are basic (Corollary~\ref{c:basic forms of the relation groupoid}).  We also show that there is a functor between $\Bi$ and $\Diffeol$ (Theorem~\ref{t:morita}), as well as that basic forms ``pullback'' by bibundles to basic forms (Proposition~\ref{p:bibundle pullback}).  These results will be needed in the sequel.

\begin{definition}[Basic Differential Forms]\labell{d:basic forms}
Given a Lie groupoid $G$, a differential form $\alpha$ on $G_0$ is \textbf{basic} if $s^*\alpha=t^*\alpha$.  Together with the standard differential, these form a complex, denoted by $\Omega^*_\basic(G)$.
\end{definition}

\begin{remark}\labell{r:PPT}
Pflaum-Posthuma-Tang define a form $\alpha\in\Omega^k(M)$ to be basic with respect to a Lie groupoid $G$ if $\rho(X)\hook\alpha=0$ for all smooth sections $X$ of the associated Lie algebroid to $G$ with anchor map $\rho$, and if $\alpha$ is $G$-invariant \cite[Definition 8.1]{PPT}.  Regarding the second condition, the first condition implies that $\alpha$ descends to a smooth section $\widetilde{\alpha}$ of the $k$th wedge power of the normal bundle to an orbit, for each orbit in $G_0$; the second condition now requires that $([s_*v]-[t_*v])\hook\widetilde{\alpha}$ vanishes for each $v\in TG_1$ (abusing notation).  It is not difficult to show that this definition is equivalent to Definition~\ref{d:basic forms} using the standard identification of the Lie algebroid with the pullback by the unit map of the subbundle $\bigcup_{x\in M}T(t^{-1}(x))\subset TG_1$; we do not need this result, however, in this paper, and so we provide no further detail.
\end{remark}

Recall that for a Lie group $K$ and a $K$-manifold $M$, a differential form on $M$ is basic if it is $K$-invariant and horizontal, the latter term meaning that it vanishes on vectors tangent to the $K$-orbits.

\begin{lemma}[Equivalence of Notions of Basicness]\labell{l:basic forms}
Let $K$ be a Lie group and $M$ a $K$-manifold.  A differential form $\alpha$ on $M$ is basic with respect to the action if and only if it is basic with respect to the action groupoid $K\ltimes M:=(K\times M\toto M)$.\footnote{The statement of this lemma was communicated to the author by Eugene Lerman; however, the proof is the author's.}
\end{lemma}

\begin{proof}
Fix $(k,x)\in K\times M$, and a vector $v\in T_{(k,x)}(K\times M)$.  Denote by $\mathfrak{k}$ the Lie algebra of $K$.  Via a left trivialisation of $TK$ we identify $T(K\times M)$ with $K\times\mathfrak{k}\times TM$.  Under this identification, there is some $\xi\in\mathfrak{k}$ and $u\in T_xM$ such that $v=(k,\xi,u)$.  It follows that $s_*v=u$ and $t_*v=\xi_M|_{k\cdot x}+k_*u$, where $\xi_M$ is the vector field on $M$ induced by $\xi\in\mathfrak{k}$.

Given $v_1,...,v_\ell\in T_{(k,x)}(K\times M)$, setting $v_i=(k,\xi_i,u_i)$ as above, for any $\ell$-form $\alpha$ on $M$:
\begin{equation}\labell{e:basic forms1}
s^*\alpha(v_1,...,v_l)=\alpha(u_1,...,u_l),
\end{equation}
\begin{equation}\labell{e:basic forms2}
t^*\alpha(v_1,...,v_l)=\alpha\big((\xi_1)_M|_{k\cdot x}+k_*u_1,...,(\xi_l)_M|_{k\cdot x}+k_*u_l\big).
\end{equation}

If follows that if $\alpha$ is $K$-invariant and horizontal, then $s^*\alpha=t^*\alpha$.

Conversely, suppose that $s^*\alpha=t^*\alpha$.  Fix $u_1,...,u_l\in T_xM$.  Since $s$ is a surjective submersion, there exist for each $i=1,...,l$ vectors $v_i=(k,\xi_i,u_i)\in K\times\mathfrak{k}\times TM$ such that $s_*v_i=u_i$.  Without loss of generality, we may take $\xi_i=0$ for each $i$.  By Equations~\eqref{e:basic forms1} and \eqref{e:basic forms2}, it follows that $\alpha$ is $K$-invariant.

Now suppose that $u_1$ is tangent to the $K$-orbit through $x$. There exists $\zeta\in\mathfrak{k}$ such that $u_1=\zeta_M|_x$.  Let $\xi_1=-\Ad_k(\zeta)$.  Then $$t_*v_1=-(\Ad_k(\zeta))_M|_{k\cdot x}+k_*u_1=0.$$  By \eqref{e:basic forms1} and \eqref{e:basic forms2} we conclude that $\alpha$ is horizontal.  The proof is complete.
\end{proof}

The following is a useful tool when checking whether a differential form on the base of a Lie groupoid is the pullback of a form on its orbit space by the quotient map.  A proof appears in \cite[Article 6.38]{iglesias}.

\begin{proposition}[Pullbacks from a Quotient]\labell{p:basic forms of the relation groupoid}
Let $X$ be a diffeological space equipped with an equivalence relation $\sim$, with quotient $Y=X/\!\sim$ and quotient map $\pi\colon X\to Y$.  For any $k$ and $\alpha\in\Omega^k(X)$, there exists $\beta\in\Omega^k(Y)$ such that $\pi^*\beta=\alpha$ if and only if $p_1^*\alpha=p_2^*\alpha$ for any plots $p_1,p_2\colon U\to X$ such that $\pi\circ p_1=\pi\circ p_2$.
\end{proposition}

\begin{remark}\labell{r:basic forms of the relation groupoid}
Proposition~\ref{p:basic forms of the relation groupoid} can be restated in terms of the relation groupoid as follows.  Let $\pr_i\colon X\times_\pi X \to X$ be the $i$th standard projection map. For any $k$ and $\alpha\in\Omega^k(X)$, there exists $\beta\in\Omega^k(Y)$ such that $\pi^*\beta=\alpha$ if and only if $\pr_1^*\alpha=\pr_2^*\alpha$.
\end{remark}

We now can prove part of Theorem~\ref{t:main}.

\begin{corollary}[Pullbacks from the Quotient are Basic]\labell{c:basic forms of the relation groupoid}
Let $G$ be a Lie groupoid.  The pullback map $\pi_G^*\colon\Omega^*(G_0/G_1)\to\Omega^*(G_0)$ is an injection with image in $\Omega_\basic^*(G)$.
\end{corollary}

\begin{proof}
Fix $\beta\in\Omega^k(G_0/G_1)$ and let $\alpha=\pi_G^*\beta$.  By Proposition~\ref{p:basic forms of the relation groupoid} and Remark~\ref{r:basic forms of the relation groupoid}, $\pr_1^*\alpha=\pr_2^*\alpha$.  Thus $(s,t)^*(\pr_1^*\alpha-\pr_2^*\alpha)=0$ (see Example~\ref{x:orbitspace}), and so $s^*\alpha=t^*\alpha$.  Injectivity follows from the definition of the quotient diffeology on $G_0/G_1$.
\end{proof}

In the case of a compact Lie group action, $\pi^*_G$ from Corollary~\ref{c:basic forms of the relation groupoid} becomes an isomorphism.  This follows from Lemma~\ref{l:basic forms} and \cite[Theorem 3.20]{watts-phd} (see also \cite{KW}), stated below.

\begin{theorem}[Group Action Case]\labell{t:KW}
Let $K$ be a compact Lie group and $M$ a $K$-manifold, with quotient map $\pi\colon M\to M/K$.  The pullback map $\pi^*$ is an isomorphism between the de Rham complexes of differential forms on $M/K$ and basic differential forms on $M$.
\end{theorem}

Sending a Lie groupoid to its diffeological orbit space constitutes a functor.

\begin{theorem}[The Functor $\Psi$]\labell{t:morita}
There is a functor $\Psi\colon\Bi\to\Diffeol$ sending a Lie groupoid $G$ to its orbit space $G_0/G_1$, a bibundle $P\colon G\to H$ to a unique smooth map $\Psi_P\colon G_0/G_1\to H_0/H_1$ such that $\Psi_P\circ\pi_G\circ a_L=\pi_H\circ a_R$, and a $2$-arrow in $\Bi$ to a trivial $2$-arrow in $\Diffeol$.
\end{theorem}

\begin{proof}
Let $P\colon G\to H$ a bibundle between Lie groupoids; denote the anchor maps to $G$ and $H$ by $a_L^P\colon P\to G_0$ (or just $a_L$ if $P$ is understood) and $a_R^P\colon P\to H_0$, respectively.
Fix $x\in G_0$.  Define $\Psi_P\colon G_0/G_1\to H_0/H_1$ by $\Psi_P([x])=\pi_H\circ a_R\circ\sigma(x)$ where $\sigma$ is a local section of $a_L$ about $x$.  To show $\Psi_P$ is well-defined, let $y$ be in the same orbit as $x$, and $\sigma'$ a local section of $a_L$ about $y$.  It suffices to show there exists $h\in H_1$ so that $a_R(\sigma'(y))=s(h)$ and $a_R(\sigma(x))=t(h)$.  Let $g\in G_1$ so that $s(g)=x$ and $t(g)=y$.  Then $a_L(g\cdot\sigma(x))=y$. Since $a_L$ is a principal $H$-bundle, there exists a unique $h$ such that $(g\cdot\sigma(x))\cdot h=\sigma'(y)$. Well-definedness of $\Psi_P$ now follows from the $G$-invariance of $a_R$.  The identity $\Psi_P\circ\pi_G\circ a_L=\pi_H\circ a_R$ follows from the definition of $\Psi_P$, and implies uniqueness.

To show that $\Psi_P$ is smooth, through local sections of $a_L$ and local lifts of a plot $p\colon U\to G_0/G_1$ (which exist by definition of the quotient diffeology), one obtains that $\Psi_P\circ p$ is locally a plot, and hence globally a plot by the Locality Axiom.

That a $2$-arrow $\alpha\colon P\Rightarrow Q$ is sent to a trivial $2$-arrow follows from the uniqueness of $\Psi_P$ and $\Psi_Q$ and the ($G$-$H$)-equivariance of $\alpha$.

Given a third Lie groupoid $K=(K_1\toto K_0)$ and bibundle $Q\colon H\to K$, the composition $Q\circ P$ is the quotient by the diagonal $H$-action $(P\times_{H_0}Q)/H$.  To show $\Psi_{Q\circ P}=\Psi_Q\circ\Psi_P$, it suffices to show $$\Psi_Q\circ\Psi_P\circ\pi_G\circ a^P_L(p) =\pi_K\circ a^Q_R(q)\quad\forall(p,q)\in P\times_{H_0}Q.$$  This follows from $\pi_H\circ a_R^P(p)=\pi_H\circ a_L^Q(q),$ which in turn follows from the definition of $P\times_{H_0}Q$. Thus $\Psi(P):=\Psi_P$ respects composition, and all other identities are straightforward to check.
\end{proof}

We now show that bibundles induce maps between complexes of basic forms, similar to what smooth maps between manifolds do with the corresponding de Rham complexes.  Versions of this result appear in the literature; for instance, \cite{PPT}, \cite{HS}.

\begin{proposition}[Pullbacks of Basic Forms by Bibundles]\labell{p:bibundle pullback}
Let $P\colon G\to H$ be a bibundle between Lie groupoids with anchor maps $a_L\colon P\to G_0$ and $a_R\colon P\to H_0$.  For any $\beta\in\Omega^k_{basic}(H)$ there exists a unique form in $\Omega^k_{basic}(G)$, denoted $P^*\beta$, such that $a_L^*(P^*\beta)=a_R^*\beta$.  In fact, $\Omega^k_\basic$ is a functor on $\Bi$ sending $P$ to a homomorphism of complexes $P^*\colon\Omega^k_\basic(H)\to\Omega^k_\basic(G)$, and a $2$-arrow in $\Bi$ to a trivial $2$-arrow.  In particular, a Morita equivalence induces an isomorphism between complexes of basic forms.
\end{proposition}

\begin{proof}

Fix $\beta\in\Omega_\basic^k(H_0)$.  The action groupoid $P\rtimes H$ has source the action map $\operatorname{act}_H$ and target $\pr_1$; with respect to these we have $$\operatorname{act}_H^*a_R^*\beta=\pr_2^*s^*\beta=\pr_2^*t^*\beta=\pr_1^*a_R^*\beta.$$ Thus $a_R^*\beta$ is basic with respect to $P\rtimes H$.

Since $P$ is a principal $H$-bundle over $G_0$, $P\rtimes H$ is isomorphic as a groupoid to $P\times_{G_0}P\toto P$.  But this is the relation groupoid for the action groupoid $G\ltimes P$, and so it follows that $a_R^*\beta$ is basic with respect to $G\ltimes P$ as well. By Remark~\ref{r:basic forms of the relation groupoid} there exists a form $\alpha$ on $G_0$ such that $a_L^*\alpha=a_R^*\beta$; $\alpha$ is unique since $a_L$ is a surjective submersion.  
Similarly, $\pr_1\colon G_1\times_{G_0}P\to G_1$ is a surjective submersion, and so the result follows from
$$\pr_1^*s^*\alpha=\pr_2^*a_L^*\alpha=\pr_2^*a_R^*\beta=\operatorname{act}_G^*a_R^*\beta=\operatorname{act}_G^*a_L^*\alpha=\pr_1^*t^*\alpha.$$

It follows from uniqueness and the definitions that $(Q\circ P)^*=P^*Q^*$ for composable bibundles $P$ and $Q$, and that isomorphic bibundles yield the same pullback map between complexes.  The remaining statements are straightforward to check.
\end{proof}

\section{Linearisations and Proof of Theorem~\ref{t:main}}\labell{s:linearisation}

In this section we review linearisations in the context of Lie groupoids, and prove Theorem~\ref{t:main}.  We follow the notation and terminology of \cite{dHF}; more details can be found in \cite[Subsection 1.2]{CS}.

Fix a Lie groupoid $G$.  A submanifold $S\subseteq G_0$ is \textbf{saturated} if it is a union of orbits of $G$.  The restriction of $G$ to $S$ is denoted by $G_S\toto S$.  The normal bundle $\nu(S)$ of $S$ in $G_0$ along with the normal bundle $\nu(G_S)$ of $G_S$ in $G_1$ form a Lie groupoid $\nu(G_S)\toto\nu(S)$ whose structure maps are induced by the differentials of those of $G$.  The standard bundle projections form a homomorphism $(\nu(G_S)\toto\nu(S))\to(G_S\toto S)$ that provides a local linear model of $G$ about $S$.

A \textbf{groupoid neighbourhood} of $G_S\toto S$ is a subgroupoid $U=(U_1\toto U_0)$ of $G$ in which $U_0$ is an open neighbourhood of $S$ and $U_1$ an open neighbourhood of $G_S$; it is \textbf{full} if $U_1=s^{-1}(U_0)\cap t^{-1}(U_0)$.  We say that $G$ is \textbf{linearisable} at a saturated submanifold $S$ if there exist full groupoid neighbourhoods $U$ of $G_S\toto S$ in $G$ and $V$ of $G_S\toto S$ in $\nu(G_S)\toto\nu(S)$, and a Lie groupoid isomorphism $\Lambda\colon V\to U$ whose restriction to $G_S$ is the identity.  A celebrated result is the following:

\begin{theorem}[\cite{CS}, Corollary 5.13 of \cite{dHF}]\labell{t:linearisable}
Proper Lie groupoids are linearisable at each of their orbits.
\end{theorem}

We now show given an orbit $\calO$ of a proper Lie groupoid $G$ that the condition that the de Rham complex of the orbit space of $\nu(G_\calO)\toto\nu(\calO)$ is isomorphic via pullback of the quotient map to basic forms of $\nu(G_\calO)\toto\nu(\calO)$ is a local one.

\begin{lemma}[Locality]\labell{l:local model forms}
Let $G$ be a proper Lie groupoid, fix $x\in G_0$, and let $\calO$ be the orbit of $x$.  If $V$ is a full groupoid neighbourhood of $\calO$ in $N:=(\nu(G_\calO)\toto\nu(\calO))$ then $\pi_V^*\colon\Omega^*(V_0/V_1)\to\Omega^*_\basic(V)$ is an isomorphism of complexes.
\end{lemma}

\begin{proof}
Let $V^N$ be the saturation of $V$ in $N$; \emph{i.e.}\ $(V^N)_0$ is the union of orbits of points in $V_0$ and $(V^N)_1=\nu(G_\calO)_{(V^N)_0}$.  By \cite[Example 3.2]{CS} $V$ and $V^N$ are Morita equivalent.  Let $G_x$ be the stabiliser of $x$ in $G$.  By \cite[Example 3.3]{CS}, there is a linear action of $G_x$ on the normal space $\nu(\calO)_x$ to $\calO$ at $x$ whose action groupoid $G_x\ltimes\nu(\calO)_x$ is Morita equivalent to $N$; this Morita equivalence restricts to one between the restriction $W$ of $G_x\ltimes\nu(\calO)_x$ to $W_0:=\nu(\calO)_x\cap(V^N)_0$ and $V^N$.  By Proposition~\ref{p:bibundle pullback} and Theorem~\ref{t:morita}, it suffices to prove $\pi_W^*\colon\Omega^*(W_0/W_1)\to\Omega^*_\basic(W)$ is an isomorphism.  But as $G_x$ is a compact Lie group, this follows from Theorem~\ref{t:KW} using $G_x$-invariant partitions of unity on $W_0$.
\end{proof}

\begin{proof}[of Theorem~\ref{t:main}]

By Corollary~\ref{c:basic forms of the relation groupoid}, it is sufficient to show that the image of $\pi_G^*$ is $\Omega^*_\basic(G)$.  Fix $\alpha\in\Omega^k_\basic(G)$.  By Proposition~\ref{p:basic forms of the relation groupoid}, it suffices to show $p_1^*\alpha= p_2^*\alpha$ for any two plots $p_1,p_2\colon U\to G_0$ such that $\pi_G\circ p_1=\pi_G\circ p_2$; fix two such plots.  It further suffices to show this equality near each point $u\in U$; fix such a point.  Let $\calO$ be the orbit of $p_1(u)$ (and hence $p_2(u)$ as well).  By Theorem~\ref{t:linearisable}, $G$ is linearisable; let $V=(V_1\toto V_0)$ and $W=(W_1\toto W_0)$ be full groupoid neighbourhoods of $G_\calO\toto\calO$ in $\nu(G_\calO)\toto\nu(\calO)$ and $G_1\toto G_0$, respectively, and let $\Lambda\colon V\to W$ be an isomorphism that fixes $G_\calO\toto\calO$.  Set $B:= p_1^{-1}(W_0)\cap p_2^{-1}(W_0)$ and let $i\colon W\to G$ be the inclusion morphism.  The pullback $\Lambda_0^*i_0^*\alpha$ is basic with respect to $V$.

By Lemma~\ref{l:local model forms}, there exists $\beta\in\Omega^k(V_0/V_1)$ such that $\pi_V^*\beta=\Lambda_0^*i_0^*\alpha$.  By Theorem~\ref{t:morita}, $\Lambda$ descends to a diffeomorphism $\Psi_\Lambda\colon V_0/V_1\to W_0/W_1$ such that $\Psi_\Lambda\circ\pi_V=\pi_W\circ\Lambda_0$, from which it follows that $\pi_W^*(\Psi^{-1}_\Lambda)^*\beta=i_0^*\alpha$.  Setting $q_n=p_n|_B$ for $n=1,2$,
$$p_n|_B^*\alpha=q_n^*i_0^*\alpha=q_n^*\pi_W^*(\Psi^{-1}_\Lambda)^*\beta.$$
The proof now is reduced to showing that $\pi_W\circ q_1=\pi_W\circ q_2$.

By Theorem~\ref{t:morita} the inclusion $i$ descends to a smooth map $j\colon W_0/W_1\to G_0/G_1$ such that for $n=1,2$,
$$j\circ\pi_W\circ q_n = \pi_G\circ i_0\circ q_n=\pi_G\circ p_n|_B.$$
Thus $j\circ\pi_W\circ q_1= j\circ\pi_W\circ q_2$, and the proof is now reduced to showing that $j$ is injective.

Fix $x_1,x_2\in W$ such that $j(\pi_W(x_1))=j(\pi_W(x_2))$.  There exists $g\in G_1$ such that $s(g)=x_1$ and $t(g)=x_2$.  Since $W$ is full, $g\in W_1$ as well, and so $\pi_W(x_1)=\pi_W(x_2)$.  This shows that $j$ is injective.

Naturality of the isomorphism of complexes follows from Proposition~\ref{p:bibundle pullback} and Theorem~\ref{t:morita}.  This completes the proof.
\end{proof}


\end{document}